\documentclass[11pt]{article}
\usepackage{latexsym}
\title{On Uniform f-vectors of Cutsets in the Truncated Boolean Lattice}

\author{B\'{e}la Bajnok\\ Gettysburg College\\ Gettysburg, PA 17325 USA \and Shahriar Shahriari\\ Pomona College\\ Claremont, CA 91711 USA}

\newtheorem{thm}{Theorem}

\newtheorem{cor}[thm]{Corollary}
\newtheorem{prop}[thm]{Proposition}
\newtheorem{conj}[thm]{Conjecture}

\newcommand{\ben}{\begin{enumerate}}
\newcommand{\een}{\end{enumerate}}
\newcommand{\bc}[2]{{{#1}\choose{#2}}}

\begin{document}

\renewcommand{\today}{August 20, 1997}

\maketitle

\begin{abstract}
Let $[n] = \{1, 2, \ldots, n\}$ and let $2^{[n]}$ be the collection of all
subsets of $[n]$ ordered by inclusion. ${\cal C} \subseteq 2^{[n]}$ is a
{\em cutset} if it meets every maximal chain in $2^{[n]}$, and the {\em width} of ${\cal C} \subseteq 2^{[n]}$ is the minimum number of chains in a chain decomposition of ${\cal C}$.  Fix $0 \leq m
\leq l \leq n$.  What is the smallest value of $k$ such that there exists
a cutset that consists only of subsets of sizes between $m$ and $l$, and
such that it contains exactly $k$ subsets of size $i$ for each $m \leq i
\leq l$?  The answer, which we denote by $g_n(m,l)$, gives a lower estimate for the width of a cutset between levels $m$ and $l$ in $2^{[n]}$.  After using the Kruskal-Katona Theorem to give a general characterization of cutsets in terms of the number and sizes of their elements, we find lower and upper bounds (as well as  some exact values) for $g_n(m,l)$.
\end{abstract}

\section{Introduction}

Let $2^{[n]}$ be the {\em Boolean lattice} of order $n$, that is the lattice of all subsets (often called {\em nodes}) of $[n] = \{1, 2, \ldots, n \}$ ordered by inclusion. For $0 \leq m \leq n$ we define the $m$-th {\em level set} $\bc{[n]}{m}$ of $2^{[n]}$ as the set of all subsets of size $m$. The {\em $f$-vector} (or {\em profile}) ${\bf f}=(f_0, f_1, \ldots, f_n)$ of a collection of subsets ${\cal A} \subseteq 2^{[n]}$ is defined by $f_m = |{\cal A}_m|$ where ${\cal A}_m = {\cal A} \cap \bc{[n]}{m}$ and $0 \leq m \leq n$.   

A collection of $l+1$ subsets $A_0 \subset A_1 \subset \cdots \subset A_l$ in $2^{[n]}$ is called a {\em chain} of length $l$. A {\em maximal chain} in $2^{[n]}$ is one that has length $n$. A collection of $w$ nodes with the property that none of them contains another is called an {\em antichain} of size $w$. The {\em length} and the {\em width} of a collection of subsets ${\cal A} \subseteq 2^{[n]}$ are defined as the length of the longest chain and the size of the largest antichain in ${\cal A}$, respectively. 

A {\em cutset} in $2^{[n]}$ is defined as a collection of subsets ${\cal C} \subseteq 2^{[n]}$ which intersects all maximal chains. Trivially, every collection ${\cal C}$ which contains $\emptyset$ or $[n]$ is a cutset. In \cite{BajnokSha:96} we proved that for $n \geq 2$, the width of a cutset in the Boolean lattice of order $n$ which does not contain $\emptyset$ or $[n]$ is greater than or equal to $n-1$, and that for $n \geq 3$ there exist cutsets of width $n-1$ in $2^{[n]}$. Thus, it is possible to construct a cutset in $2^{[n]}$ with  $f$-vector $(0, \underbrace{n-1, n-1, \ldots, n-1}_{n-1}, 0)$. We then may ask for the smallest value of $k$ for which there is a cutset in $2^{[n]}$ with $f$-vector $(0, \underbrace{k, k, \ldots, k}_{n-1}, 0)$. The original goal of our work was to show that this value is $n-2$ (see Corollary \ref{cor} below). 

More generally, for $0 \leq m \leq l \leq n$ we define $g_n(m,l)$ to be the smallest value of $k$ for which the $n+1$-tuple $(f_0, f_1, \ldots, f_n)$, defined by $f_i=k$ if $m \leq i \leq l$ and 0 otherwise, can be the $f$-vector of a cutset in $2^{[n]}$. Thus our goal above is then to find $g_n(1,n-1)$. Note that by symmetry we have $g_n(m,l)=g_n(n-l,n-m)$, so we may assume without loss of generality that $m \leq l \leq n-m$. 

Before studying $g_n(m,l)$, we give a general characterization of $f$-vectors of cutsets in $2^{[n]}$. For a given profile ${\bf f}=(f_0,f_1,\dots,f_n)$ and integer $m_0$ with $0 \leq m_0 \leq n$, we construct a {\em canonical} collection of subsets ${\cal C}({\bf f},m_0)$, with the property that there is a cutset in $2^{[n]}$ with profile ${\bf f}$ if and only if ${\cal C}({\bf f},m_0)$ is a cutset for some (or every) $0 \leq m_0 \leq n$. We then translate this qualitative criterion to a quantitative one: For a given ${\bf f}=(f_0,f_1,\dots,f_n)$ and $0 \leq m_0 \leq n$, we describe an easily computable value $q({\bf f},m_0)$, so that ${\bf f}$ will be the profile of a cutset in $2^{[n]}$ exactly when $f_{m_0} \geq q({\bf f},m_0)$ for some (or every) $0 \leq m_0 \leq n$. These characterizations, which we present in Section 2, are essentially due to Daykin \cite{Daykin:75b} (for a correction see \cite{Clements:84} and then \cite{Clements:?}), though we follow a treatment which is more suitable for our purposes.

We can then determine the values of $g_n(m,l)$ for $l \leq m+2$. Namely, we prove the following.   

\begin{thm} \label{values}

Let $n$ be a positive integer.

\ben

	\item

$g_n(m,m)=\bc{n}{m}$ for every integer $0 \leq m \leq n$.

	\item

$g_n(m,m+1)=\bc{n-1}{m}$ for every integer $0 \leq m \leq n-1$.

	\item

$g_n(m,m+2)=\sum_{j=0}^{m} \bc{n-2j-2}{m-j}$  for every integer $0 \leq m \leq n/2-1$.

\een

\end{thm} 

Next, viewing $m$ as fixed and $n>>m$ (i.e., for all $n>n_0=n_0(m)$), we develop upper and lower bounds for $g_n(m,l)$.

\begin{thm} \label{bounds}

Suppose that $m$ and $n$ are non-negative integers and $n >> m$. Then

\begin{enumerate}
	\item

$\bc{n-2}{m} < g_n(m,l) \leq \sum_{j=0}^{m} \bc{n-2j-2}{m-j}$ for every integer $m+2 \leq l \leq n-m-1$, and

	\item

$\bc{n-3}{m} < g_n(m,n-m) \leq \sum_{j=0}^{m} \bc{n-2j-3}{m-j}$.

\end{enumerate}

\end{thm}

For $m=1$ we then get the following results.

\begin{cor} \label{cor}

Suppose that $n>4$ and $1 \leq l \leq n-1$ are integers. Then
$$g_n(1,l)=\left\{ \begin{array}{ll} n & \mbox{if $l=1$} \\
				     n-1 & \mbox{if $2 \leq l \leq n-2$} \\
				     n-2 & \mbox{if $l=n-1$}
		\end{array}
		\right.$$
\end{cor}

For $2 \leq m << n$, Theorems \ref{values} and \ref{bounds} give the ``numerator'' of the leading term of the $m$-binomial representation (see section 2) of $g_n(m,l)$. Namely, this value is equal to $n$ if $l=m$, $n-1$ if $l=m+1$, $n-2$ if $m+2 \leq l \leq n-m-1$, and $n-3$ if $l=n-m$. It is striking that for a rather large range of values of $l$, $g_n(m,l)$ stays essentially unchanged.

We note that in Theorem \ref{bounds}, the ratio of the upper bound to the lower bound is approximately $1+\frac{m}{n}$, and thus the bounds are rather accurate as $n>>m$.

	Extremal problems regarding cutsets in the Boolean lattice have been the object of much study.  For example see \cite{DuffusSanWin:90,FurGriKle89,Griggs:95,GriKle89,Kleitman:94,Now87,ShiLi93}.

\section{$f$-vectors of cutsets}

Given a collection ${\cal B} \subseteq \bc{[n]}{m}$, the {\em shadow} and the {\em shade} of ${\cal B}$ will be denoted by $\bigtriangleup {\cal B}$
and $\bigtriangledown {\cal B}$, respectively \cite[Chapter 2]{And89}, and are as usual  defined by $$ \bigtriangleup {\cal B} = \{ A \in \bc{[n]}{m-1} \mid A \subseteq B \; \ \mbox{for some} \; \ B \in {\cal B}\},$$ $$ \bigtriangledown {\cal B} = \{ A \in \bc{[n]}{m+1} \mid B \subseteq A \; \ \mbox{for
some} \; \ B \in {\cal B}\}.$$        

We order the elements of $\bc{[n]}{m}$ by the {\em squashed order} (also called the {\em colex} order) \cite[Chapter 7]{And89}, that is for $A, B \in \bc{[n]}{m}$, we say $A <_S B$ if the largest element of the symmetric difference of $A$ and $B$ is in $B$. For $1 \leq K \leq \bc{n}{m}$, we define the {\em initial collection} ${\cal F}_m(K)$ and the {\em last collection} ${\cal L}_m(K)$ at level $m$ as the first and last $K$ elements in the squashed order at level $m$, respectively. In addition, if $K \leq 0$, then ${\cal F}_m(K)={\cal L}_m(K)=\emptyset$. The squashed order has the property that the shadow of an initial collection at level $m$ is an initial collection at level $m-1$, and the shade of a last collection at level $m$ is a last collection at level $m+1$. The Kruskal-Katona Theorem (\cite{Kat66b,Kru63} or \cite[Chapter 7]{And89}) states that the size of the shadow of $K$ nodes at level $m$ is greater than or equal to the size of the shadow of ${\cal F}_m(K)$ and, equivalently, the size of their shade is greater than or equal to the size of the shade of ${\cal L}_m(K)$.

Let $\Omega_n$ denote the set of $n+1$-tuples of integers $(a_0,a_1,\dots,a_n)$ such that $0 \leq a_m \leq \bc{n}{m}$ for all $0 \leq m \leq n$. To see whether a given ${\bf f} \in \Omega_n$ is the profile of a cutset in $2^{[n]}$, we construct a collection of subsets ${\cal C} = {\cal C}({\bf f},m_0)$, called the {\em canonical} collection of subsets for profile ${\bf f}$ and for level $m_0$ ($0 \leq m_0 \leq n$). As we show below,  there is a cutset in $2^{[n]}$ with profile ${\bf f}$ if and only if this canonical collection is a cutset for some (or every) $m_0$. 

Our construction is as follows. First we let ${\cal E}_{0}^{\uparrow} = \{ \emptyset \}$, ${\cal C}_{0}^{\uparrow} = {\cal F}_0(f_0)$, and for $1 \leq m \leq n$ we recursively define ${\cal E}_{m}^{\uparrow}  = \bigtriangledown ({\cal E}_{m-1}^{\uparrow} - {\cal C}_{m-1}^{\uparrow})$ and ${\cal C}_{m}^{\uparrow}  = {\cal L}_m(|{\cal E}_{m}^{\uparrow}|) - {\cal L}_m(|{\cal E}_{m}^{\uparrow}| - f_m)$. Then ${\cal E}_{m}^{\uparrow}$ is a last collection at level $m$, and it is precisely the set of nodes from which there is a chain of length $m$ to $\emptyset$ which is disjoint from ${\cal C}_{i}^{\uparrow}$ for all $0 \leq i \leq m-1$. Analogously, we let ${\cal E}_{n}^{\downarrow} = \{[n]\}$, ${\cal C}_{n}^{\downarrow} = {\cal L}_n(f_n)$, and for $0 \leq m \leq n-1$ we recursively define ${\cal E}_{m}^{\downarrow}  = \bigtriangleup ({\cal E}_{m+1}^{\downarrow} - {\cal C}_{m+1}^{\downarrow})$ and ${\cal C}_{m}^{\downarrow}  = {\cal F}_m(|{\cal E}_{m}^{\downarrow}|) - {\cal F}_m(|{\cal E}_{m}^{\downarrow}| - f_m)$. This time ${\cal E}_{m}^{\downarrow}$ is an initial collection at level $m$, and it is the set of nodes from which there is a chain of length $n-m$ to $[n]$ which is disjoint from ${\cal C}_{i}^{\uparrow}$ for all $m+1 \leq i \leq n$. Finally, we define ${\cal C} = {\cal C}({\bf f},m_0) = (\cup_{i=0}^{m_0} {\cal C}_{i}^{\uparrow}) \cup (\cup_{i=m_0+1}^{n}{\cal C}_{i}^{\downarrow}$).

We can easily see that ${\cal C}={\cal C}({\bf f},m_0)$ is a cutset if and only if ${\cal E}_{m_0}^{\uparrow} \cap {\cal E}_{m_0}^{\downarrow} \subseteq {\cal C}_{m_0}$. Furthermore, the profile $(c_0, c_1, \ldots, c_n)$ of ${\cal C}$ satisfies $c_m \leq f_m$ for every $m$, and if ${\cal C}$ is not a cutset, then its profile is exactly ${\bf f}$.

For example, let $n = 5$.  If  ${\bf f} = (0,2,5,6,0,0)$ then, for all $0 \leq m_0 \leq
5$, ${\cal C}({\bf f},m_0)$ becomes $$\{
\{1\},\{2\},\{1,3\},\{2,3\},\{1,4\},\{2,4\},\{3,4\},\{1,2,5\},\{1,3,5\},\{2,3,5\},\{1,4,5\},\{2,4,5\},\{3,4,5\}\}.$$ 
On the other, hand if  ${\bf g} = (0,2,6,5,0,0)$, then for ${\cal C}({\bf g},5)$ we get
$$\{
\{1\},\{2\},\{1,3\},\{2,3\},\{1,4\},\{2,4\},\{3,4\},\{1,5\},\{1,2,5\},\{1,3,5\},\{2,3,5\},\{1,4,5\},\{2,4,5\}\}.$$
It is easily seen that the first is a cutset, while the second one is a
collection with profile ${\bf g}$ and not a cutset.

The next two propositions give us useful ways of determining whether a given vector ${\bf f} \in \Omega_n$ can be the profile of a cutset in $2^{[n]}$.

\begin{prop} \label{canonical}

Let ${\bf f} \in \Omega_n$ and $0 \leq m_0 \leq n$. The canonical collection ${\cal C}({\bf f},m_0)$ defined above is a cutset if and only if $|{\cal E}_{m_0}^{\uparrow}|+|{\cal E}_{m_0}^{\downarrow}| \leq \bc{n}{m_0}+f_{m_0}$.
\end{prop} 

{\em Proof.\/} If ${\cal C}$ is a cutset, then the assertion follows as $|{\cal E}_{m_0}^{\uparrow}|+|{\cal E}_{m_0}^{\downarrow}| = 
|{\cal E}_{m_0}^{\uparrow} \cap {\cal E}_{m_0}^{\downarrow}| + |{\cal E}_{m_0}^{\uparrow} \cup {\cal E}_{m_0}^{\downarrow}|
\leq |{\cal C}_{m_0}| + \bc{n}{m_0} \leq f_{m_0} + \bc{n}{m_0}$.

On the other hand, if ${\cal C}$ is not a cutset, then its profile must be exactly ${\bf f}$. Furthermore, as ${\cal E}_{m_0}^{\downarrow}$ is an initial segment and ${\cal E}_{m_0}^{\uparrow}$ is a last segment at level $m_0$, their intersection has size greater than $f_{m_0}$ and their union is all of $\bc{[n]}{m_0}$. Therefore, in this case we have $|{\cal E}_{m_0}^{\uparrow}|+|{\cal E}_{m_0}^{\downarrow}| = 
|{\cal E}_{m_0}^{\uparrow} \cap {\cal E}_{m_0}^{\downarrow}| + |{\cal E}_{m_0}^{\uparrow} \cup {\cal E}_{m_0}^{\downarrow}|
> f_{m_0} + \bc{n}{m_0}$. $\qquad \Box$

The next proposition shows the importance of the canonical collection. It further shows that the choice of $m_0$ is immaterial.

\begin{prop} \label{qualitative}

Let ${\bf f} \in \Omega_n$. The following are equivalent.
\ben \item\ ${\bf f}$ is the profile of a cutset in $2^{[n]}$.
\item\   ${\cal C}({\bf f},m_0)$ is a cutset for some $0 \leq m_0 \leq n$.
\item\   ${\cal C}({\bf f},m_0)$ is a cutset for every $0 \leq m_0 \leq n$.\een
\end{prop} 

{\em Proof.\/}  Clearly, if ${\cal C}({\bf f},m_0)$ is a cutset for some $0 \leq m_0 \leq n$, then  there is a cutset in $2^{[n]}$ with profile ${\bf f}$. Therefore, it is enough to prove that if there is an ${\cal A} \subseteq 2^{[n]}$ which is a cutset with profile ${\bf f}$, then ${\cal C}({\bf f},m_0)$ is a cutset for every $0 \leq m_0 \leq n$.

For $m = 0, \ldots n$, let  ${\cal A}_m = \{ A \in {\cal A} \mid |A| = m
\}$. We define ${\cal B}_{m}^{\uparrow}$ and ${\cal B}_{m}^{\downarrow}$ recursively by ${\cal B}_{0}^{\uparrow} = \{ \emptyset \}, {\cal B}_{m}^{\uparrow}  =
\bigtriangledown ({\cal B}_{m-1}^{\uparrow} - {\cal A}_{m-1})$ for $m = 1, 2, \ldots, n$, and 
 ${\cal B}_{n}^{\downarrow} = \{ [n] \}, {\cal B}_{m}^{\downarrow}  = \bigtriangleup ({\cal B}_{m+1}^{\downarrow} - {\cal A}_{m+1})$ for $m = 0,1, \ldots, n-1$. Since ${\cal A}$ is a cutset, we must have ${\cal B}_{m}^{\downarrow} \cap {\cal B}_{m}^{\uparrow} \subseteq {\cal A}_m$ for every $0 \leq m \leq n$. 

Keeping the notations established prior to Proposition \ref{canonical}, we now use downward induction on $m$ to prove that $|{\cal B}_{m}^{\downarrow}| \geq |{\cal E}_{m}^{\downarrow}|$ and $|{\cal B}_{m}^{\uparrow}| \geq |{\cal E}_{m}^{\uparrow}|$ for every $0 \leq m \leq n$. Clearly $|{\cal B}_{n}^{\downarrow}| \geq |{\cal E}_{n}^{\downarrow}|$. Using, in order, the definition of ${\cal B}_{m}^{\downarrow}$, the Kruskal-Katona theorem, the triangle inequality, the inductive hypothesis, the fact that ${\cal C}_{m+1}^{\downarrow} \subseteq {\cal E}_{m+1}^{\downarrow}$, that ${\cal E}_{m+1}^{\downarrow} - {\cal C}_{m+1}^{\downarrow}$ is an initial segment at level $m+1$, and the definition of ${\cal E}_{m}^{\downarrow}$, we can write 
\begin{eqnarray*}
|{\cal B}_{m}^{\downarrow}| & = &|\bigtriangleup ({\cal B}_{m+1}^{\downarrow} - {\cal A}_{m+1})| \\
	& \geq & |\bigtriangleup {\cal F}_{m+1}(|{\cal B}_{m+1}^{\downarrow} - {\cal A}_{m+1}|)| \\
	& \geq & |\bigtriangleup {\cal F}_{m+1}(|{\cal B}_{m+1}^{\downarrow}| - |{\cal A}_{m+1}|)| \\
	& \geq & |\bigtriangleup {\cal F}_{m+1}(|{\cal E}_{m+1}^{\downarrow}| - |{\cal C}_{m+1}^{\downarrow}|)| \\
	& = & |\bigtriangleup {\cal F}_{m+1}(|{\cal E}_{m+1}^{\downarrow} - {\cal C}_{m+1}^{\downarrow}|)| \\
	& = & |\bigtriangleup ({\cal E}_{m+1}^{\downarrow} - {\cal C}_{m+1}^{\downarrow})| \\
	& = & |{\cal E}_{m}^{\downarrow}|,
\end{eqnarray*}
 as claimed. The assertion for $|{\cal B}_{m}^{\uparrow}|$ can be proved similarly.

Our assertion now follows from Proposition \ref{canonical}, as 

$|{\cal E}_{m_0}^{\uparrow}| + |{\cal E}_{m_0}^{\downarrow}| \leq |{\cal B}_{m_0}^{\uparrow}| + |{\cal B}_{m_0}^{\downarrow}| = |{\cal B}_{m_0}^{\uparrow} \cap {\cal B}_{m_0}^{\downarrow}| + |{\cal B}_{m_0}^{\uparrow} \cup {\cal B}_{m_0}^{\downarrow}| \leq |{\cal A}_{m_0}| + |{\cal B}_{m_0}^{\uparrow} \cup {\cal B}_{m_0}^{\downarrow}| \leq  f_{m_0} + \bc{n}{m_0}.$ $\qquad \Box$

A quantitative version of Proposition \ref{qualitative} can be formulated as follows.

Given positive integers $K$ and $m$, there exist unique integers $a_m > a_{m-1} > \cdots > a_t \geq t \geq 1$ such that 
$$K = \bc{a_m}{m} + \bc{a_{m-1}}{m-1} + \cdots + \bc{a_t}{t}.$$
This is called the {\em $m$-binomial representation} of $K$ \cite[Theorem 7.2.1]{And89}. Using the $m$-binomial representation of $K$ it is easy to describe the set numbered $K$ in the squashed order on the $m$-th level of the Boolean lattice \cite[page 117]{And89}.   

 For each positive integer $m$ we define a (boundary) operator $\partial_m$ \cite{GreKle78} (\cite{Col87} and \cite{FranklMatRuzTok:95} have other notations) on the integers as follows: If $K$ is a positive integer with an $m$-binomial representation as above then
$$\partial_m(K) =  \bc{a_m}{m-1} + \bc{a_{m-1}}{m-2} + \cdots + \bc{a_t}{t-1},$$ 
and for non-positive $K$ set $\partial_m(K) =0$. Note that $\partial_m$ is weakly increasing.

With this operator we can write $\bigtriangleup{\cal F}_m(K) = {\cal F}_{m-1}(\partial_m(K))$, and $\bigtriangledown{\cal L}_m(K) = {\cal L}_{m+1}(\partial_{n-m}(K))$, and the Kruskal-Katona theorem becomes $|\bigtriangleup{\cal B}| \geq \partial_m(|{\cal B}|)$ and $|\bigtriangledown{\cal B}| \geq \partial_{n-m}(|{\cal B}|)$ where {\cal B} is a collection of subsets of $[n]$ at level $m$.

\begin{thm} \label{quantitative}

For a given ${\bf f}=(f_0, f_1, \ldots, f_n) \in \Omega_n$, define ${\bf u(f)}=(u_0, u_1, \ldots, u_n)$ and ${\bf v(f)}=(v_0, v_1, \ldots, v_n)$ by:

\hspace{1in} $u_0 = 1$ and $u_{m+1} = \partial_{n-m}(u_{m}-f_{m})$ for $m = 0, \ldots, n-1,$

\hspace{1in} $v_n = 1$ and $v_{m-1} = \partial_{m}(v_{m}-f_{m})$ for $m = 1, \ldots, n.$

Then the following are equivalent.
\ben \item\ ${\bf f}$ is the profile of a cutset in $2^{[n]}$.
\item\   $u_m + v_m - f_m \leq \bc{n}{m}$  for some $0 \leq m \leq n$.
\item\   $u_m + v_m - f_m \leq \bc{n}{m}$  for every $0 \leq m \leq n$.
\een
\end{thm}

{\em Proof.\/} $(u_0, u_1, \ldots, u_n)$ and $(v_0, v_1, \ldots, v_n)$ were chosen so that $|{\cal E}_{m}^{\uparrow}| = u_m$ and $|{\cal E}_{m}^{\downarrow}| = v_m$, so the statements follow from Propositions \ref{canonical} and \ref{qualitative}. $\qquad \Box$

Note that, as a special case when $f_0=f_n=0$, ${\bf f}$ is the profile of a cutset in $2^{[n]}$ if and only if $v_0=0$ (or if and only if $u_n=0$).

{\em Remark.\/} We can give a quantitative description of the canonical collection ${\cal C}({\bf f},m_0) = (\cup_{i=0}^{m_0} {\cal C}_{i}^{\uparrow}) \cup (\cup_{i=m_0+1}^{n}{\cal C}_{i}^{\downarrow}$) as follows. If levels of $2^{[n]}$ are in the squashed order, then ${\cal C}_{i}^{\uparrow}$ is the segment at level $i$ starting with node numbered $\bc{n}{i} - u_i+1$ and ending with $\min\{\bc{n}{i} - u_i+f_i,\bc{n}{i}\}$; ${\cal C}_{i}^{\downarrow}$ is the segment at level $i$ starting with node numbered $\max\{1,v_i-f_i+1\}$ and ending with $v_i$.

\section{Exact values: Proof of Theorem \ref{values}}

We will use Theorem \ref{quantitative} to determine the values of $g_n(m,l)$ for $l=m$, $m+1$, and $m+2$. 

The case $l=m$ is obvious. For $l=m+1$, we can easily see that the canonical collection with profile ${\bf f}=(0,0, \dots, 0, \bc{n-1}{m},\bc{n-1}{m}, 0, \dots,0)$ (where the nonzero components occur at levels $m$ and $m+1$) consists of subsets of size $m$ that do not contain $n$, together with subsets of size $m+1$ that contain $n$. This collection is a cutset. On the other hand, the canonical collection with profile ${\bf f'}=(0,0, \dots, 0, \bc{n-1}{m},\bc{n-1}{m}-1, 0, \dots,0)$ is not a cutset as there exists a maximal chain through the node $\{n-m,n-m+1,\dots,n\}$ which does not intersect the collection.
	
For the case $l = m+2$, the result is trivial for $m = 0$ and hence assume
that $m \geq 1$.  Define $f = \sum_{j=0}^{m-1} \bc{n-2j-2}{m-j}$, and we
need to prove that $g_n(m,l) = f+1$.  We will show that ${\bf f}=(0,0,
\dots, 0, f+1,f,f, 0, \dots,0)$ (where the nonzero components occur at
levels $m$, $m+1$, and $m+2$) is the profile of a cutset, but that ${\bf
f'}=(0,0, \dots, 0, f,f,f, 0, \dots,0)$ is not the profile of a cutset.
 
Before starting the proof we need to establish a binomial identity which establishes a relationship between two \emph{vertical} columns in the arithmetic (a.k.a. Pascal's) triangle.  Let
$n$ be a positive integer,  let $1 \leq m \leq n/2$, and let $0 \leq d
\leq m-1$.  Using downward induction on $d$ (with base case $d=m-1$) it can be easily shown that
$$\sum_{j=0}^{m-d-1} \bc{n-2j-2}{m-j}+\sum_{j=0}^{m-d-1}
\bc{n-2j-1}{m+2-j}+\bc{n-2m+2d}{d+2}=\bc{n}{m+2}.$$
In our computations below we rely on the case $d = 0$ of this identity, namely:
$$f+\sum_{j=0}^{m-1} \bc{n-2j-1}{m+2-j}+\bc{n-2m}{2}=\bc{n}{m+2}.$$ 

We now use Theorem \ref{quantitative} to show that ${\bf f}$ is the profile of a cutset. Computing ${\bf v(f)}$ we get

\begin{eqnarray*}
v_i & = & \bc{n}{i} \qquad \mbox{ for } i=n, n-1,\dots, m+2, \\
v_{m+1} & = & \partial_{m+2}[\bc{n}{m+2}-f]\\
& = & \partial_{m+2}[\sum_{j=0}^{m-1} \bc{n-2j-1}{m+2-j}+\bc{n-2m}{2}] \\
 & = & \sum_{j=0}^{m-1} \bc{n-2j-1}{m+1-j}+\bc{n-2m}{1}, \\
v_{m} & = & \partial_{m+1}[\sum_{j=0}^{m-1} \bc{n-2j-1}{m+1-j}+\bc{n-2m}{1}-f]\\
&= & \partial_{m+1}[\sum_{j=0}^{m-1} \bc{n-2j-2}{m+1-j}+\bc{n-2m}{1}] \\
 & = & \sum_{j=0}^{m-1} \bc{n-2j-2}{m-j}+\bc{n-2m}{0}\\
& = & f+1, \\
v_{m-1} & = & \partial_{m}[f+1-f-1]\\
& = & 0, \qquad \mbox{ and } \\
v_i & = & 0 \qquad \mbox{ for } i=m-1,m-2,\dots,0, 
\end{eqnarray*}
hence ${\bf f}$ is the profile of a cutset. As ${\bf f'}$ differs from ${\bf f}$ only in its $m$-th coordinate, we still have $v_m=f+1$. Continuing, we get

\begin{eqnarray*}
v_{m-1} & = & \partial_{m}[f+1-f]\\
& = & \partial_{m}[\bc{m}{m}]\\
& = & \bc{m}{m-1}, \mbox{ and } \\
v_i & = & \bc{m}{i} \qquad \mbox{ for } i=m-1,m-2,\dots,0. 
\end{eqnarray*}

In particular, $v_0=1$, so ${\bf f'}$ is not the profile of a cutset. $\qquad \Box$

{\em Remark.\/} We can also construct a (non-canonical) cutset with $f$-vector ${\bf f}=(0,0, \dots, 0, f+1,f,f, 0, \dots,0)$ where $f = \sum_{j=0}^{m-1} \bc{n-2j-2}{m-j}$. We define ${\cal Q}_0=\bc{[n-2]}{m}$ and for $j=1,2,\dots,m$ we let 
\begin{eqnarray*}
{\cal Q}_j & = & \{A \cup \{n-1,n-3,\dots,n-2j+1\}|A \in \bc{[n-2j-2]}{m-j} \},\\
{\cal R}_j & =& \{Q \cup \{n-2j+2\} | Q \in {\cal Q}_{j-1} \}, \mbox{and} \\
{\cal S}_j &= & \{R \cup \{n-2j+1\} | R \in {\cal R}_{j-1} \}.
\end{eqnarray*}

Then it is rather straight forward to check that ${\cal Q}= \bigcup_{j=0}^{m}{\cal Q}_j$, ${\cal R}= \bigcup_{j=1}^{m}{\cal R}_j$, and ${\cal S}= \bigcup_{j=1}^{m}{\cal S}_j$ are collections of nodes at levels $m$, $m+1$, and $m+2$, respectively, $|{\cal Q}|=f+1$, $|{\cal R}|=f$, and $|{\cal S}|=f$, and that ${\cal Q} \cup {\cal R} \cup {\cal S}$ is a cutset, thus providing a cutset with $f$-vector ${\bf f}=(0,0, \dots, 0, f+1,f,f, 0, \dots,0)$. We omit the details but only briefly sketch the idea behind the
construction:  Note that ${\cal Q}_0$ consists of {\em all} subsets of
size $m$ that do not contain $n$ or $n-1$.  Likewise, ${\cal R}_1$
consists of all subsets of size $m+1$ that contain $n$ but not $n-1$ and
${\cal S}_1$ consists of all subsets of size $m+2$ that contain $n$ and
$n-1$.  Thus any maximal chain in $2^{[n]}$ that does not intersect ${\cal
Q}_0 \cup  {\cal R}_1 \cup {\cal S}_1$ will have subsets of size $m$, $m+1$ and
$m+2$ that contain $n-1$ but not $n$.  Now the poset of subsets of $[n]$
that contain $n-1$ but not $n$ is isomorphic (as a poset) to $2^{[n-2]}$.
We now restrict ourself to this poset and continue the construction
recursively.

\section{Bounds: Proof of Theorem \ref{bounds}}

We need the following propositions.

\begin{prop} \label{complement}

For $m \leq l < n/2$ we have $g_n(m,n-m) \leq g_{n-1}(m,l)$.

\end{prop}

{\em Proof.\/} Let ${\cal A} \subseteq 2^{[n-1]}$ be a cutset with $f$-vector $(0,\dots,0,g,\dots,g,0,\dots,0)$, where $g=g_{n-1}(m,l)$ and the nonzero entries are at levels $m \leq i \leq l$. Let ${\cal \overline{A}} = \{[n] \setminus A | A \in {\cal A} \}$.  Note that by symmetry ${\cal \overline{A}}$ is a cutset in the poset ${\cal Q}$ consisting of all subsets of $[n]$ that contain $n$.
  Now define
 ${\cal B}={\cal A}  \cup \overline{{\cal A}}$.
 Note that since $l<n-l$, ${\cal B}$ has $f$-vector $(0,\dots,0,g,\dots,g,0,\dots,0,g,\dots,g,0,\dots,0)$, where the nonzero entries are at levels $m \leq i \leq l$ and $n-l \leq i \leq n-m$.
It now suffices to show that ${\cal B}$ is a cutset in $2^{[n]}$.

Let ${\cal C}=C_0 \subset C_1 \subset \cdots \subset C_n$ be a maximal chain in $2^{[n]}$ with $|C_i|=i$ for $0 \leq i \leq n$, and suppose, indirectly, that ${\cal B} \cap {\cal C} = \emptyset$. Since ${\cal A}$ is a cutset in $2^{[n-1]}$ with support between levels $m$ and $l$, we must have $n \in C_k$ for every $l \leq k \leq n$.  Similarly, since ${\cal \overline{A}}$ is a cutset in ${\cal Q}$, we must have $n \not\in C_k$ for every $0 \leq k \leq n-l$.  This can only happen if $n-l < l$ which is a contradiction. $\qquad \Box$

\begin{prop} \label{doubling}

Suppose that $m+2 \leq l \leq n-m-1$ and $g_{n-1}(m+1,l-1)>\bc{n-2}{m}$. Then  $g_{n}(m,l)>\bc{n-2}{m}$.

\end{prop}

{\em Proof.\/} Let ${\bf f}=(0,0,\dots,0,\bc{n-2}{m},\bc{n-2}{m},\dots,\bc{n-2}{m},0,\dots,0)$, where the nonzero entries occur between levels $m$ and $l$. Suppose, indirectly, that ${\bf f}$ is the profile of a cutset in $2^{[n]}$. Then, by Proposition \ref{qualitative}, the canonical collection ${\cal C}({\bf f},n)$ is such a cutset. If $C_i= {\cal C} \cap \bc{[n]}{i}$ are the level sets of the collection, then it is easy to see that $C_m = \bc{[n-2]}{m}$ and $C_{m+1} = \{A \cup \{n-1\} | A \in \bc{[n-2]}{m}\}$. Therefore $\cup_{i=m+2}^{l}C_i$ is a (canonical) cutset in $\{A \cup \{n\} | A \subseteq [n-1]\}$ of profile $(0,0,\dots,0,\bc{n-2}{m},\bc{n-2}{m},\dots,\bc{n-2}{m},0,\dots,0)$, where the nonzero entries occur between levels $m+2$ and $l$. But $\{A \cup \{n\} | A \subseteq [n-1]\}$ is isomorphic to $2^{[n-1]}$, so we must have $g_{n-1}(m+1,l-1) \leq \bc{n-2}{m}$, a contradiction. $\qquad \Box$

\begin{prop} \label{minusing}

Define ${\bf f}=(f_0,f_1,\dots,f_n)$ as follows. Let $f_{m}=f_{m+1}=f_{n-m-1}=f_{n-m}=\bc{n-3}{m}$, $f_i=f_{n-i}=\bc{n+m-i-1}{m+1}$ for $i=m+2,m+3,\dots,\lfloor n/2 \rfloor$, and $f_i=0$ otherwise. Then ${\bf f}$ is not the profile of a cutset in $2^{[n]}$.

\end{prop}

{\em Proof.\/} For ${\bf v(f)}$ we get the following.

\begin{eqnarray*}
v_i & = & \bc{n}{i} \qquad \mbox{for} \quad i=n,n-1,\dots,n-m, \\
v_{n-m-1} & = & \partial_{n-m}[\bc{n}{n-m}-\bc{n-3}{n-m-3}]\\
	&=& \partial_{n-m}[\bc{n-1}{n-m}+\bc{n-2}{n-m-1}+\bc{n-3}{n-m-2}] \\
 & = & \bc{n-1}{n-m-1}+\bc{n-2}{n-m-2}+\bc{n-3}{n-m-3}, \\
v_{n-m-2} & = & \partial_{n-m-1}[\bc{n-1}{n-m-1}+\bc{n-2}{n-m-2}+\bc{n-3}{n-m-3}-\bc{n-3}{n-m-3}] \\
	& = & \partial_{n-m-1}[\bc{n-1}{n-m-1}+\bc{n-2}{n-m-2}]\\
& = & \bc{n-1}{n-m-2}+\bc{n-2}{n-m-3}, \\
v_{n-m-3}& = & \partial_{n-m-2}[\bc{n-1}{n-m-2}+\bc{n-2}{n-m-3}-\bc{n-3}{n-m-4}] \\
	& = & \partial_{n-m-2}[\bc{n-1}{n-m-2}+\bc{n-3}{n-m-3}]\\
& = & \bc{n-1}{n-m-3}+\bc{n-3}{n-m-4}, \\
v_{n-m-4} & = & \partial_{n-m-3}[\bc{n-1}{n-m-3}+\bc{n-3}{n-m-4}-\bc{n-4}{n-m-5}] \\
	& = & \partial_{n-m-3}[\bc{n-1}{n-m-3}+\bc{n-4}{n-m-4}]\\
& = & \bc{n-1}{n-m-4}+\bc{n-4}{n-m-5}, 
\end{eqnarray*}
and continuing as above, we get $v_k=\bc{n-1}{k}+\bc{k+m}{k-1}$ for every $k \geq \lfloor(n-1)/2 \rfloor$.

Now note that since ${\bf f}$ is symmetrical, for ${\bf v(f)}$ and ${\bf u(f)}$ we have $v_i=u_{n-i}$ for every $0 \leq i \leq n$. Therefore

$v_{\lfloor n/2 \rfloor}+u_{\lfloor n/2 \rfloor}
=v_{\lfloor n/2 \rfloor}+v_{\lceil n/2 \rceil}
=\bc{n-1}{\lfloor n/2 \rfloor}+\bc{\lfloor n/2 \rfloor+m}{\lfloor n/2 \rfloor-1}
+\bc{n-1}{\lceil n/2 \rceil}+\bc{\lceil n/2 \rceil+m}{\lceil n/2 \rceil-1}
=\bc{n}{\lfloor n/2 \rfloor}+f_{\lfloor n/2 \rfloor}+f_{\lceil n/2 \rceil}$,
thus ${\bf f}$ cannot be the profile of a cutset by Theorem \ref{quantitative}. $\qquad \Box$

{\em Proof of Theorem \ref{bounds} and Corollary \ref{cor}.\/} The upper bound in 1 follows from Theorem \ref{values} part 3 since $g_n(m,l)$ is weakly decreasing with $l$. The upper bound in 2 is a consequence of Proposition \ref{complement} (take $l=m+2$ assuming $n>2m+4$) and Theorem \ref{values} part 3. Proposition \ref{minusing} implies that $g_n(m,n-m)>\min\{\bc{n-3}{m},\bc{\lfloor n/2 \rfloor +m-1}{m+1}\}$. It is an elementary exercise to check that this minimum is $\bc{n-3}{m}$ for $n>>m$ (e.g. if $n \geq 3^{m+1}(m+1)$), establishing the lower bound in 2. Finally, the lower bound in 1 is a consequence of Proposition \ref{doubling} taking $l=n-m-1$, and noting that $\bc{n-4}{m+1} \geq \bc{n-2}{m}$ for $n>>m$ (e.g. if $n \geq 8(m+1)$) and the lower bound in 2. 

The cases of $l=1,2$ (and $3$) of Corollary \ref{cor} follow from Theorem \ref{values}. The assertions for $3 \leq l \leq n-1$ follow from Theorem \ref{bounds} if $n \geq 18$ (see above). The cases $ 5 \leq n \leq 17$ were checked directly using Theorem \ref{quantitative} (and a simple computer program). 
$\qquad \Box$

We close this section by proving a partial complement to Proposition \ref{complement}.

\begin{prop} \label{symm}

For $n/2-1 < l \leq n-m-2$ we have $g_n(m,n-m) \geq g_{n-1}(m,l)$.

\end{prop}

{\em Proof.\/} Consider the canonical collection ${\cal C}={\cal C}({\bf f},l)$ where ${\bf f}=(f_0,f_1,\dots,f_n)$ is defined by $f_i=g_n(m,n-m)$ for $m \leq i \leq n-m$ and 0 otherwise. Then ${\cal C}$ is a cutset in $2^{[n]}$. Write ${\cal D}_1=\{C \in {\cal C} | |C| \leq l\}$ and ${\cal D}_2=\{C \in {\cal C} | |C| \geq l+1\}$. 

Suppose, indirectly, that $g_n(m,n-m) < g_{n-1}(m,l)$. Then ${\cal D}_1$ is not a cutset in $2^{[n-1]}$, in particular, ${\cal D}_1 \subseteq 2^{[n-1]}$. Since $l \geq n-1-l$, by symmetry ${\cal D}_2$ is disjoint from $2^{[n-1]}$. Therefore ${\cal D}_1 \cup {\cal D}_2 = {\cal C}$ is not a cutset in $2^{[n]}$, a contradiction. $\qquad \Box$

\section{Examples and conjectures}

The following table has the exact values of $g_n(m,l)$ for $n = 100$, $m =
4$ and every $4 \leq l \leq 96$.  These were found using Theorem 6.  We
have also given the $4$-binomial representation of $g_n(m,l)$ in the third
column.

$$\begin{array}{||c||c||c||} \hline \hline
l & g_{100}(4,l) & \mbox{4-binomial representation of $g_{100}(4,l)$} \\
\hline \hline
4 &3,921,225 & \bc{100}{4} \\ \hline
5 &3,764,376 & \bc{99}{4} \\ \hline
6 & 3,759,624& \bc{98}{4}+\bc{96}{3}+\bc{94}{2}+\bc{93}{1}\\ \hline
7 &3,759,526 & \bc{98}{4}+\bc{96}{3}+\bc{93}{2}+\bc{88}{1} \\ \hline
8 \leq l \leq 95 & 3,759,525 &
\bc{98}{4}+\bc{96}{3}+\bc{93}{2}+\bc{87}{1}=\bc{100}{4}-\bc{100}{3} \\
\hline
96 & 3,607,527 & \bc{97}{4}+\bc{95}{3}+\bc{92}{2}+\bc{86}{1}=
\bc{99}{4}-\bc{99}{3}\\ \hline \hline
\end{array}$$
Theorem 1 gives the exact value for $4 \leq l \leq 6$.  For $6 \leq l \leq
95$, Theorem 2 gives $3,612,280 = \bc{98}{4} < g_{100}(4,l) \leq
\bc{98}{4}+\bc{96}{3}+\bc{94}{2}+\bc{92}{1}+\bc{90}{0}=3,759,624$.
Finally, for $l = 96$, Theorem 2 gives $3,464,840 =\bc{97}{4} <
g_{100}(4,96) \leq
\bc{97}{4}+\bc{95}{3}+\bc{93}{2}+\bc{91}{1}+\bc{89}{0}=3,607,625$.

From these and other similar tables we see how $g_n(m,l)$ decreases as $l$
increases from $m$ to $n-m$. Namely, we observe that the decrease is
largest from level $m$ to level $m+1$ and from level $n-m-1$ to level
$n-m$, quite modest between level $m+2$ and level $2m$, and that, rather
strikingly, $g_n(m,l)$ is constant between levels $l=2m$ and $l=n-m-1$. In
fact, we have the following conjectures.

\begin{conj} \label{conj}

For $n>>m$ we have

\begin{enumerate}

	\item

$g_n(m,l)=\bc{n}{m}-\bc{n}{m-1}$ for every $l=2m,2m+1, \dots,n-m-1$, and

\item

$g_n(m,n-m)=\bc{n-1}{m}-\bc{n-1}{m-1}$.

\end{enumerate}

\end{conj}

Note that the $m$-binomial representation of $\bc{n}{m}-\bc{n}{m-1}$ starts with $\bc{n-2}{m}+\bc{n-4}{m-1}$ (and $\bc{n-1}{m}-\bc{n-1}{m-1}$ starts with $\bc{n-3}{m}+\bc{n-5}{m-1}$) when $n>>m$ (cf. Theorem \ref{bounds}). Corollary \ref{cor} proves our conjectures for $m=1$, and Theorem \ref{bounds} establishes that Conjecture \ref{conj} provides an upper bound for $m=2$ since $\bc{n}{2}-\bc{n}{1}=\bc{n-2}{2}+\bc{n-3}{1}$ (and $\bc{n-1}{2}-\bc{n-1}{1}=\bc{n-3}{2}+\bc{n-4}{1}$).

According to Conjecture \ref{conj}, we have $g_n(m,n-m)=g_{n-1}(m,n-m-2)=g_{n-1}(m,l)$  for $n>>m$ and $2m \leq l \leq n-m-2$. Proposition \ref{complement} establishes $g_n(m,n-m) \leq g_{n-1}(m,l)$ for $l<n/2$, while Proposition \ref{symm} proves the other direction for $n/2-1 < l \leq n-m-2$ (yielding equality when $l=(n-1)/2$).

{\em Note:\/} In a subsequent paper \cite{Bajnok:98}, the above
conjectures have been somewhat refined and related to other conjectures
about the width of cutsets in the truncated Boolean lattice.

\end{document}